Theo van den Bogaart, august 2008

# About the choice of a basis in Kedlaya's algorithm

The goal of this text is to provide a proof of Proposition 5.3.1 in [Edix], which is a statement concerning a part of Kedlaya's algorithm. This algorithm will not be explained here; it is described in [Ked] and a detailed introduction can be found in the course notes [Edix]. Both these references exclude characteristic 2. See [D-V] for an adaptation of the algorithm to this situation.

Let $p$ be a prime number. Kedlaya's algorithm calculates the number of points of an arbitrary hyperelliptic curve over a finite field of characteristic $p$. It does this by determining the Zeta function of the curve. The key idea is that by the theory of Monsky and Washnitzer, the Frobenius morphism of the curve induces an automorphism of the de Rham cohomology of a lift of the curve over some $p$-adic field. The usual cohomological machinery then produces the Zeta function. The algorithm calculates a $p$-adic approximation of the Frobenius, in the form of a matrix with respect to a specific basis for the de Rham cohomology. In this text, the denominators appearing in the coefficients of this matrix are investigated.

More precisely, start with a hyperelliptic curve $\bar{C}$ over a finite field $\mathbb{F}_q$ with $q = p^n$ elements. In Kedlaya's algorithm, one chooses a lift of $\bar{C}$ to a curve $C$ over the ring of Witt vectors $\mathbb{Z}_q$ of $\mathbb{F}_q$. The Frobenius endomorphism of $\bar{C}$ canonically induces a map on the de Rham cohomology $\mathrm{H}^1_{\mathrm{dR}}(C'_{\mathbb{Q}_q}/\mathbb{Q}_q)$, where $C' \subset C$ is a certain open affine subscheme, where $C'_{\mathbb{Q}_q} = C' \otimes \mathbb{Q}_q$ and where $\mathbb{Q}_q$ is the quotient field of $\mathbb{Z}_q$. This action stabilises the $(-1)$-eigenspace for the hyperelliptic involution $\mathrm{H}^1_{\mathrm{dR}}(C'_{\mathbb{Q}_q}/\mathbb{Q}_q)^-$. It is this last space that is used in Kedlaya's algorithm. A basis $B$ is specified and the matrix $M$ corresponding to the $p$-power Frobenius is calculated by the algorithm. One then calculates, semi-linearly, the $n$th power of this matrix to obtain the $q$-Frobenius on de Rham cohomology. In general, the coefficients of $M$ will not lie in $\mathbb{Z}_q$.

In section 2 it will be shown that $\mathrm{H}^1_{\mathrm{dR}}(C/\mathbb{Z}_q)$ is a lattice in $\mathrm{H}^1_{\mathrm{dR}}(C'_{\mathbb{Q}_q}/\mathbb{Q}_q)^-$ that is Frobenius-stable. This essentially follows from the comparison theorem with crystalline cohomology. Hence changing the basis from $B$ to a basis $B'$ that generates the $\mathbb{Z}_q$-module $\mathrm{H}^1_{\mathrm{dR}}(C/\mathbb{Z}_q)$ results in a matrix $M'$ which has $\mathbb{Z}_q$-coefficients. The problem how to make this transition is treated in sections 3 and 4. The key to this is the gap sequence for hyperelliptic curves.

In section 5 all this is then applied to Kedlaya's algorithm, assuming the characteristic $p$ to be different from 2. It is shown that $p^{[\log_p(2g-1)]}M$ has integral coefficients, and an algorithm is exhibited that passes from $B$ to $B'$. Finally, in section 6, there are some comments on the case $p=2$; the results obtained in this situation are less satisfactory.

*Remark:* Instead of $\mathrm{H}^1_{\mathrm{dR}}(C/\mathbb{Z}_q) \subset \mathrm{H}^1_{\mathrm{dR}}(C'_{\mathbb{Q}_q}/\mathbb{Q}_q)$ one can also consider the image of

$$\mathrm{H}^1_{\mathrm{dR}}(C'/\mathbb{Z}_q) \longrightarrow \mathrm{H}^1_{\mathrm{dR}}(C'_{\mathbb{Q}_q}/\mathbb{Q}_q)$$

and it seems to me an interesting question whether or not this gives a Frobenius stable lattice. See also the third remark on page 205 in [M-W].

Because we need to keep track of the integral structure on $C$, we cannot afford the luxury to reduce everything to explicit calculations with geometric objects, i.e., curves defined over a field. In particular, we need to work with 'families of curves' over $\mathbb{Z}_q$ and we need some technical statements about how the cohomology of such objects compares to the cohomology of the fibres. In this spirit, section 1 recalls some theory about families of curves over an arbitrary base. In the following sections, generality gradually decreases: while section 2 is valid for higher dimensional objects too, section 3 is particular for curves and in section 4 we specialize further to the case of hyperelliptic curves. Sections 5 and 6 depend on explicit calculations with the affine equation of a hyperelliptic curve.

To calculate the number of fixed points of a curve over a field with $p^n$ elements, one has to calculate the characteristic polynomial of the $p^n$-power Frobenius. Therefore, one has to take the $n$th power of the matrix of Frobenius (in a semi-linear way and upto a sufficiently large power of $p$). Kedlaya's algorithm is polynomial in the genus of the curve and the degree of the finite field over $\mathbb{F}_p$, but not polynomial in $\log p$. So in practice, one often takes $p$ small, but $n$ very large. Therefore, with respect to memory costs it is good to have a matrix that is denominator-free.

Finally some words about the history of the problem at hand. That the basis used in Kedlaya's original algorithm can give a matrix with non-integral coefficients was noticed by Vercauteren when he implemented the algorithm. The main theorem of this text has appeared in the course notes [Edix] of Edixhoven about Kedlaya's algorithm. In these notes there is only a rough sketch of the proof. Kedlaya refers to Edixhoven's theorem in [Ked$^\flat$]. Edixhoven asked the author to work out a precise proof based on his ideas, which resulted in this text.

*Notation:* Fix an integer $q = p^n$, with $p$ a prime and $n \geq 1$. Note that $p = 2$ is allowed. Let $k$ be a finite field with $q$ elements, $\mathcal{V} = W(k)$ its ring of Witt vectors and let $K = \mathcal{V}[\frac{1}{p}]$ be the quotient field of $\mathcal{V}$. (For aesthetic reasons we have decided not to adopt the notations $k = \mathbb{F}_q$, $\mathcal{V} = \mathbb{Z}_q$ and $K = \mathbb{Q}_q$ used in the introduction and in [Ked].) Let $\sigma$ be the Teichmüller lift to $\mathcal{V}$ and $K$ of the $p$-power Frobenius automorphism $x \mapsto x^p$ of $k$.

If $C$ is a scheme over $\mathcal{V}$, the notation $C_k = C \times_{\operatorname{Spec} \mathcal{V}} \operatorname{Spec} k$ is used for the special fibre and $C_K = C \times_{\operatorname{Spec} \mathcal{V}} \operatorname{Spec} K$ for the generic fibre. The *Frobenius morphism* of a scheme will always mean the $p$-power Frobenius (so $x \mapsto x^p$ on the coordinates).

We fix an integer $g \geq 1$, which plays the rôle of genus.

Let the reader beware that the notation used in this text is slightly different from the one used in [Ked]. Most importantly: the letter $C$ is used in [Ked] to denote a curve over $k$; here it is used for a lift to $\mathcal{V}$.

## 1. Preliminaries on hyperelliptic curves

This text is involved with hyperelliptic curves defined over $\mathcal{V}$. Such curves are obtained for example by lifting to $\mathcal{V}$ the defining equations for a hyperelliptic curve over $k$, as in Kedlaya's algorithm. As the properties of hyperelliptic curves over a more arbitrary base than a field are perhaps less well-known, in this section we will recall some facts about families of hyperelliptic curves over an arbitrary locally noetherian base scheme $S$. (For the purpose of this text, the reader can safely specialize to $S = \operatorname{Spec} \mathcal{V}$ or the spectrum of a field.) The reference for all this, at least if $g \geq 2$, is [L-K].

A *curve over* $S$ is a smooth, projective morphism of schemes $C \to S$ whose geometric fibres are connected of dimension 1. We will also say that $C/S$ is a curve. If all the geometric fibres have genus $g'$ for some $g' \in \mathbb{Z}$, then $C/S$ is called a *curve of genus* $g'$.

Let $C/S$ be a curve of genus $g$ (recall that $g \geq 1$ by assumption). An $S$-morphism $\iota\colon C \to C$ is a *hyperelliptic involution* if $\iota^2 = \mathrm{id}$ and if for each geometric point $\bar{\eta} \to S$ there is an isomorphism $C_{\bar{\eta}}/\langle \iota_{\bar{\eta}} \rangle \simeq \mathbb{P}^1_{k(\bar{\eta})}$, where $\iota_{\bar{\eta}}$ is the induced involution on the fibre $C_{\bar{\eta}}$. It follows (Th. 4.12 in [ibid.]) that $D = C/\langle \iota \rangle$ is a curve of genus 0 (in fact $D_{\bar{\eta}} = C_{\bar{\eta}}/\langle \iota_{\bar{\eta}} \rangle$) and that the projection $C \to D$ is surjective, flat and finite of degree 2. The map $C_{\bar{\eta}} \to D_{\bar{\eta}}$ gives a separable extension of the corresponding function fields. If $g \geq 2$ the existence of a surjective, finite map of degree 2 from $C$ to a curve of genus 0 is equivalent to the existence of a hyperelliptic involution; and this involution is unique ([ibid.], Th. 5.5). A curve of genus $g$ is called a *hyperelliptic curve of genus* $g$ if there exists a hyperelliptic involution.

Let $C/S$ be a hyperelliptic curve of genus $g$. Assume (for simplicity) that $S$ is reduced. A section $s{:}S \to C$ corresponds to a relative effective Cartier divisor $\mathcal{D}$ on $C/S$ of degree 1 (see [K-M], Ch. 1). Such a section is a *Weierstrass point* if $\iota s = s$ for a hyperelliptic involution $\iota$, or equivalently, if for every geometric point $\bar{\eta} \to S$ the $k(\bar{\eta})$-vector space $\mathrm{H}^0(C_{\bar{\eta}}, \mathcal{O}_{C_{\bar{\eta}}}(2\mathcal{D}_{\bar{\eta}}))$ has dimension 2, where $\mathcal{D}_{\bar{\eta}}$ is the pull-back of $\mathcal{D}$ to the fibre above $\bar{\eta}$. A Weierstrass point $s$ determines a point of $D = C/\langle\iota\rangle$, with corresponding relative divisor $\mathcal{D}'$ on $D$. For such a Weierstrass point, the canonical maps

$$\mathrm{H}^0(\mathbb{P}^1_{k(\bar{\eta})}, \mathcal{O}_{\mathbb{P}^1_{k(\bar{\eta})}}(n\mathcal{D}'_{\bar{\eta}})) \longrightarrow \mathrm{H}^0(C_{\bar{\eta}}, \mathcal{O}_{C_{\bar{\eta}}}(2n\mathcal{D}_{\bar{\eta}})) \longrightarrow \mathrm{H}^0(C_{\bar{\eta}}, \mathcal{O}_{C_{\bar{\eta}}}((2n+1)\mathcal{D}_{\bar{\eta}}))$$

are isomorphisms for $n = 0, 1, \dots, g-1$ and for every geometric point $\bar{\eta} \to S$. (This behaviour is encoded in the *gap sequence* for hyperelliptic curves: $1, 3, \dots, 2g-1$. This is the list of positive integers $r$ for which there is no function that is regular everywhere except at one Weierstrass point where it has a pole of exact order $r$.)

## 2. Lifting Frobenius to de Rham cohomology

First we recall the definition of de Rham cohomology from [Gr] or from §7 of [Ha70]. Starting with a smooth morphism of schemes $f{:}X \to S$, let

$$\text{(2.1)} \qquad \Omega^\bullet_{X/S}: \quad 0 \longrightarrow \mathcal{O}_X \longrightarrow \Omega^1_{X/S} \longrightarrow \Omega^2_{X/S} \longrightarrow \cdots$$

be its de Rham complex (with $\mathcal{O}_X$ in degree 0). Note that in general $\Omega^\bullet_{X/S}$ is not a complex in the category of $\mathcal{O}_X$-modules, as the maps need only be $\mathcal{O}_S$-linear. If $X$ is affine, de Rham cohomology is just given by the closed $S$-differential forms modulo the exact ones, i.e., as the homology of the complex

$$\Gamma(X, \Omega^\bullet_{X/S}): \quad 0 \longrightarrow \Gamma(X, \mathcal{O}_X) \longrightarrow \Gamma(X, \Omega^1_{X/S}) \longrightarrow \Gamma(X, \Omega^2_{X/S}) \longrightarrow \cdots.$$

In general, one uses hypercohomology: the (absolute) de Rham cohomology of $X/S$ is the hypercohomology of the de Rham complex: $\mathrm{H}^i_{\mathrm{dR}}(X/S) := \mathrm{H}^i(X, \Omega^\bullet_{X/S})$ for $i \in \mathbb{Z}$. (A quick introduction to hypercohomology can be found in Appendix C of [Mi]; a detailed treatment is EGA $0_{\mathrm{III}}$ §11. There is also an enlightening Čech construction. We will not use it, but it can be used to give alternative proofs to some of the statements in this text. See for instance [K-O], page 205–206.) One can filter the de Rham complex by replacing $\Omega^i_{X/S}$ by 0 for $i$ small. The filtered complex induces a spectral sequence

$$E_1^{a,b} = \mathrm{H}^b(C, \Omega^a_{X/S}) \qquad \Rightarrow \qquad \mathrm{H}^{a+b}_{\mathrm{dR}}(X/S),$$

which in this case is called the *Hodge to de Rham spectral sequence.* If $X$ is affine, this spectral sequence degenerates since all higher order cohomology groups vanish (EGA III.1.3.1; or [HAG] III.3.7 in the noetherian case), thus showing that the two descriptions of de Rham cohomology in the affine case are compatible. We will need the following base change property.

**2.2. Proposition.** *Suppose $C$ is a smooth, separated scheme of finite type over $\mathcal{V}$. For each $i \in \mathbb{Z}$, the natural map*

$$\mathrm{H}^i_{\mathrm{dR}}(C/\mathcal{V}) \otimes_{\mathcal{V}} K \longrightarrow \mathrm{H}^i_{\mathrm{dR}}(C_K/K)$$

*is an isomorphism.*

PROOF. The relevant Hodge-de Rham spectral sequences are

$$E_1^{a,b} = \mathrm{H}^b(C, \Omega^a_{C/\mathcal{V}}) \qquad \Rightarrow \qquad \mathrm{H}^{a+b}_{\mathrm{dR}}(C/\mathcal{V})$$

and

$$'E_1^{a,b} = \mathrm{H}^b(C_K, \Omega^a_{C_K/K}) \qquad \Rightarrow \qquad \mathrm{H}^{a+b}_{\mathrm{dR}}(C_K/K).$$

There is a canonical map, obtained by pull-back along the projection $X_K \to X$, from the first spectral sequence to the last. As $K$ is flat over $\mathcal{V}$, we can tensor the first spectral sequence with $K$ and obtain a spectral sequence that abuts to $\mathrm{H}^{p+q}_{\mathrm{dR}}(C/\mathcal{V}) \otimes_{\mathcal{V}} K$. By flat base change ([HAG] Prop. 9.3 or EGA III.1.4.15) the induced map $E_1^{a,b} \otimes \mathbb{Q} \to {}'E_1^{a,b}$ is an isomorphism. But a map of spectral sequences that is an isomorphism on $E_1$ is an isomorphism of spectral sequences. □

Suppose $C$ is a proper and smooth scheme over $\mathcal{V}$. The de Rham cohomology of $C$ only depends on the reduction $C_k$ of $C$ to $k$. In fact, for each $i \in \mathbb{Z}$ there is a canonical isomorphism

$$\mathrm{H}^i_{\mathrm{dR}}(C/\mathcal{V}) \xrightarrow{\sim} \mathrm{H}^i_{\mathrm{cris}}(C_k/\mathcal{V})$$

with the crystalline cohomology of $C_k$; see [Ber1]. By functoriality, the Frobenius morphism of $C_k$ induces a $\sigma$-linear automorphism of $\mathrm{H}^i_{\mathrm{dR}}(C/\mathcal{V})$. Tensoring with $K$ we extend this automorphism to $\mathrm{H}^i_{\mathrm{dR}}(C_K/K)$.

In the situation of Kedlaya's algorithm, one considers an affine open subscheme $C' \subset C$. By making a comparison with Monsky-Washnitzer cohomology, the Frobenius morphism also induces a map on $\mathrm{H}^i_{\mathrm{dR}}(C'_K/K)$. We want to know if the restriction map

$$\mathrm{H}^i_{\mathrm{dR}}(C_K/K) \longrightarrow \mathrm{H}^i_{\mathrm{dR}}(C'_K/K) \tag{2.3}$$

is compatible with the maps induced by Frobenius on both sides. To answer this question, we need a cohomology theory for $k$-varieties which compares to crystalline cohomology in the proper case, to Monsky-Washnitzer cohomology in the affine case and which has the right functorial behaviour. This is provided by Berthelot's rigid cohomology. We will now make a small detour to this rigid cohomology in order to prove (see the theorem at the end of this section) that (2.3) is Frobenius equivariant under some reasonable assumptions.

A good introduction to rigid spaces and their cohomology is [F-vdP], especially the sections 7.6 and 7.7. For the precise definitions and main properties of rigid cohomology, we refer to [Ber2]. The relation between rigid cohomology and ordinary de Rham cohomology is explained in [B-C]. We will now describe the maps we need and how they fit together.

Let $C' \subset C$ be an arbitrary open subscheme of the proper and smooth $\mathcal{V}$-scheme $C$. Associated to $C_K$ and $C'_K$ are rigid analytic spaces $C_K^{\mathrm{an}}$ and $C_K'^{\mathrm{an}}$, respectively. In order to state the definition of rigid cohomology, we need the formal completion $\hat{C} = (C_k, \varprojlim(\mathcal{O}_C/p^n\mathcal{O}_C))$ of $C$ along $C_k$. The rigid analytic space $\hat{C}_K$ (the 'generic fibre in the sense of Raynaud') associated to it is in this situation just $\hat{C}_K = C_K^{\mathrm{an}}$ (which also equals $]C_k[$). It comes with a map $\mathrm{sp}\colon \hat{\mathrm{C}}_{\mathrm{K}} \to \hat{\mathrm{C}}$ of ringed spaces. Let $]C'_k[ := \mathrm{sp}^{-1}(C'_k) \subset \hat{C}_K$ be the tube of $C'_k$. There are inclusions

$$]C'_k[ \subset C_K'^{\mathrm{an}} \subset C_K^{\mathrm{an}} = \hat{C}_K.$$

Let $\Omega^\bullet_{C_K'^{\mathrm{an}}}$ be the analytic de Rham complex on $C_K'^{\mathrm{an}}$ (which is defined similarly to (2.1)). Define on $\hat{C}_K$ the complex

$$\Omega^{\bullet\dagger}_{C'} = \varinjlim_{]C'_k[\subset V} (j_V)_* \Omega^\bullet_{C_K'^{\mathrm{an}}}|_V,$$

where the limit is taken over all subspaces $V \subset C_K'^{\mathrm{an}}$ that are strict neighbourhoods of $]C'_k[$ in $\hat{C}_K$ and where $j_V\colon V \hookrightarrow \hat{C}_K$ is the inclusion.
The various inclusions induce morphisms

$$\begin{array}{ccccc} \mathrm{H}^i(C_K^{\mathrm{an}}, \Omega^\bullet_{C_K^{\mathrm{an}}}) & \longrightarrow & \mathrm{H}^i(C_K'^{\mathrm{an}}, \Omega^\bullet_{C_K'^{\mathrm{an}}}) & \longrightarrow & \mathrm{H}^i(\hat{C}_K, \Omega^{\bullet\dagger}_{C'}) \\ \|{\scriptstyle (\mathrm{def})} & & & & \|{\scriptstyle (\mathrm{def})} \\ \mathrm{H}^i_{\mathrm{rig}}(C_k/K) & & & & \mathrm{H}^i_{\mathrm{rig}}(C'_k/K) \end{array}.$$

The spaces on the left and right only depend on the reductions $C_k$ and $C'_k$ of $C$ and $C'$, and are by definition the rigid cohomologies of $C_k$ and $C'_k$, respectively. The map of ringed spaces $C_K'^{\mathrm{an}} \to C'_K$ induces a map $\mathrm{H}^i_{\mathrm{dR}}(C'_K/K) \to \mathrm{H}^i(C_K'^{\mathrm{an}}, \Omega^\bullet_{C_K'^{\mathrm{an}}})$. Likewise, we obtain a map $\mathrm{H}^i_{\mathrm{dR}}(C_K/K) \to \mathrm{H}^i(C_K^{\mathrm{an}}, \Omega^\bullet_{C_K^{\mathrm{an}}})$, which is an isomorphism by the non-archimedean GAGA-theorem (it also follows from the general theorem below). These two maps form a commutative diagram with the restriction maps:

$$
\begin{array}{ccc}
\mathrm{H}^i_{\mathrm{dR}}(C_K/K) & \xrightarrow{\sim} & \mathrm{H}^i_{\mathrm{rig}}(C_k/K) \\
\downarrow & & \downarrow \\
\mathrm{H}^i_{\mathrm{dR}}(C'_K/K) & \longrightarrow & \mathrm{H}^i_{\mathrm{rig}}(C'_k/K)
\end{array} \quad . \tag{2.4}
$$

By the functorial properties of rigid cohomology, the map in the right side of this square commutes with Frobenius. The next theorem tells us that sometimes the bottom map in (2.4) is also an isomorphism.

**2.5. Theorem (Kiehl, Baldassarri-Chiarellotto [B-C], Cor. 2.6).** *Suppose $C \setminus C'$ is given by a relative divisor with normal crossings. Then the map*

$$\mathrm{H}^i_{\mathrm{dR}}(C'_K/K) \longrightarrow \mathrm{H}^i_{\mathrm{rig}}(C'_k/K)$$

*is an isomorphism for all $i \in \mathbb{Z}$.* □

Note that the rigid cohomology of a smooth affine $k$-scheme is canonically isomorphic to the Monsky-Washnitzer cohomology of that scheme (see [Ber2], Prop. 1.10). Likewise, there is a natural isomorphism ([ibid.], Prop. 1.9)

$$\mathrm{H}^i_{\mathrm{rig}}(C_k/K) \xrightarrow{\sim} \mathrm{H}^i_{\mathrm{cris}}(C_k/\mathcal{V}) \otimes_{\mathcal{V}} K.$$

We conclude:

**2.6. Theorem.** *Let $C/\mathcal{V}$ be a smooth and proper scheme and let $C' \subset C$ be an open affine subscheme such that its complement is a relative divisor with normal crossings. Let $i \in \mathbb{Z}$. Then $\mathrm{H}^i_{\mathrm{dR}}(C'_K/K)$ is canonically isomorphic to the Monsky-Washnitzer cohomology of $C'_k$. If we use this to let the Frobenius of $C'_k$ act on $\mathrm{H}^i_{\mathrm{dR}}(C'_K/K)$, and if we use the comparison with crystalline cohomology to get a Frobenius on $\mathrm{H}^i_{\mathrm{dR}}(C/\mathcal{V})$, then the restriction map*

$$\mathrm{H}^i_{\mathrm{dR}}(C/\mathcal{V}) \longrightarrow \mathrm{H}^i_{\mathrm{dR}}(C'_K/K).$$

*is Frobenius equivariant.* □

## 3. Finding a lattice

Throughout this section, we fix a curve $C$ over $\mathcal{V}$ of genus $g$, together with an effective relative divisor $\mathcal{D}$ of degree 1, i.e., a divisor corresponding to a $\mathcal{V}$-valued point.

Roughly speaking, in Kedlaya's algorithm the basis chosen for de Rham cohomology (with rational coefficients) is given by differentials on $C$ having poles in $\mathscr{D}$ only, of degree $\leq 2g$. The purpose of this section is to get a criterium how far off such differentials are in forming a set of generators for $\mathrm{H}^1_{\mathrm{dR}}(C/\mathscr{V})$. This will result in a technical result, Lemma 3.10 below. The reader is encouraged to read ahead the statements in paragraph 4 to see where we are heading.

First we need two results about the cohomology of curves over $\mathscr{V}$:

**3.1. Lemma.** *Let $i \in \mathbb{Z}$ and let $\mathscr{D}'$ be a relative divisor on $C$. The $\mathscr{V}$-modules $\mathrm{H}^i(C, \mathscr{O}_C)$ and $\mathrm{H}^i(C, \Omega^1_{C/\mathscr{V}})$ are free. If $\deg \mathscr{D}' > 0$ then $\mathrm{H}^i(C, \Omega^1_{C/V}(\mathscr{D}'))$ is free. If $\deg \mathscr{D}' > 2g-2$ then $\mathrm{H}^i(C, \mathscr{O}_C(\mathscr{D}'))$ is free.*

Proof. Let $\mathscr{F} = \mathscr{O}_C(\mathscr{D}')$ or $\mathscr{F} = \Omega_{C/\mathscr{V}}(\mathscr{D}')$ (for the first statement of the lemma we can take $\mathscr{D}' = 0$). By the criterion for cohomological flatness (EGA III.7.8.4) it suffices to show that the dimension of $\mathrm{H}^i(C_s, \mathscr{F}_s)$ is independent of $s \in \operatorname{Spec} \mathscr{V}$. The Euler characteristic $\dim \mathrm{H}^0(C_s, \mathscr{F}_s) - \dim \mathrm{H}^1(C_s, \mathscr{F}_s)$ is independent of $s$ (EGA III.7.9.4), so it suffices to show that one of the two terms appearing in the Euler characteristic is independent of $s$. But for $s \in \operatorname{Spec} \mathscr{V}$ we have $\dim \mathrm{H}^0(C_s, \mathscr{O}_{C_s}) = 1$ and $\dim \mathrm{H}^0(C_s, \Omega^1_{C_s/k(s)}) = g$, which concludes the first statement. Now we use Serre duality: if $\deg \mathscr{D}' > 0$ then $\dim \mathrm{H}^1(C_s, \Omega^1_{C_s/k(s)}(\mathscr{D}'_s)) = 0$, and if $\deg \mathscr{D}' > 2g-2$ then $\dim \mathrm{H}^1(C_s, \mathscr{O}_{C_s}(\mathscr{D}'_s)) = 0$, which gives the last two statements. □

**3.2. Proposition.** *The $\mathscr{V}$-module $\mathrm{H}^1_{\mathrm{dR}}(C/\mathscr{V})$ is a lattice in $\mathrm{H}^1_{\mathrm{dR}}(C_K/K)$.*

Proof. By Proposition 2.2, it suffices to show that $\mathrm{H}^1_{\mathrm{dR}}(C/\mathscr{V})$ is free. Consider the 'stupid' short exact sequence

$$0 \longrightarrow \Omega^1_{C/\mathscr{V}}[1] \longrightarrow \Omega^\bullet_{C/\mathscr{V}} \longrightarrow \mathscr{O}_C[0] \longrightarrow 0.$$

(For an object $\mathscr{F}$ and $n \in \mathbb{Z}$, the notation $\mathscr{F}[n]$ stands for the complex which has $\mathscr{F}$ in degree $n$ and is zero elsewhere.) A part of the induced long exact sequence reads

$$\mathrm{H}^0(C, \mathscr{O}_C) \longrightarrow \mathrm{H}^0(C, \Omega^1_{C/\mathscr{V}}) \longrightarrow \mathrm{H}^1_{\mathrm{dR}}(C/\mathscr{V}) \longrightarrow \mathrm{H}^1(C, \mathscr{O}_C).$$

The space on the left is free and hence, as $C_K$ is complete, it is isomorphic to $\mathscr{V}$. The left-most map in this sequence is the derivation, which equals the zero map as $d\mathscr{V} = 0$. Now if $x \in \mathrm{H}^1_{\mathrm{dR}}(C/\mathscr{V})$ were torsion, it must map to 0 in $\mathrm{H}^1(C, \mathscr{O}_C)$, as this space is free by the lemma above. Hence $x \in \mathrm{H}^0(C, \Omega^1_{C/\mathscr{V}}) \subset \mathrm{H}^1_{\mathrm{dR}}(C/\mathscr{V})$ is torsion, which implies $x = 0$ as this space is also free. □

Consider the map

$$\mathscr{O}_C((2g-1)\mathscr{D}) \xrightarrow{d} \Omega^1_{C/\mathscr{V}}(2g\mathscr{D}) \tag{3.3}$$

obtained by restriction of the derivation $K(C) \to \Omega^1_{K(C)/K}$. Compose it with the canonical projection to the quotient sheaf $\Omega^1_{C/\mathscr{V}}(2g\mathscr{D}) \,/\, \Omega^1_{C/\mathscr{V}}(\mathscr{D})$. The cokernel of the resulting map

$$\mathscr{O}_C\big((2g-1)\mathscr{D}\big) \longrightarrow \frac{\Omega^1_{C/\mathscr{V}}(2g\mathscr{D})}{\Omega^1_{C/\mathscr{V}}(\mathscr{D})} \tag{3.4}$$

is a skyscraper sheaf whose support is contained in the support of the divisor $\mathscr{D}$.

**3.5. Notation.** The $\mathscr{V}$-module that corresponds to the cokernel of (3.4) is denoted by $\Upsilon = \Upsilon_{C,\mathscr{D}}$. There is a natural map

$$\mathrm{H}^0(C, \Omega^1_{C/\mathscr{V}}(2g\mathscr{D})) \xrightarrow{\varphi} \Upsilon \tag{3.6}$$

denoted by $\varphi = \varphi_{C,\mathscr{D}}$.

The importance of the map $\varphi$ is that its *kernel* is related to the de Rham cohomology of $C$. The relation between de Rham cohomology and this kernel will be made precise in the lemma below. But first we give a very concrete description of $\varphi$.

For this, let $t$ be a generator for the ideal $\mathscr{O}_C(-\mathscr{D}) \subset \mathscr{O}_C$ in a neighbourhood of $\operatorname{Supp}\mathscr{D}$ and let $P$ be the special point of $\operatorname{Supp}\mathscr{D}$. Then the completion $\mathscr{O}^{\wedge}_{C,P}$ of the local ring of $C$ at $P$ with respect to its maximal ideal is $\mathscr{V}[\![t]\!]$ (see Bourbaki, AC VIII.5, Th. 2). Taking the stalk of (3.4) in $P$ and completing, we obtain the derivative map

$$t^{-(2g-1)}\mathscr{V}[\![t]\!] \xrightarrow{d} \bigoplus_{i=-2g}^{-2} \mathscr{V}\, t^i dt. \tag{3.7}$$

Now localisation and completion are exact, and the codomain of this map is killed by a power of the maximal ideal in $P$. So the cokernel of (3.7) is $\Upsilon$ and we compute

$$\Upsilon = \bigoplus_{i=-2g}^{-2} \big(\mathscr{V}/(i+1)\mathscr{V}\big)\, t^i dt \simeq \bigoplus_{\substack{-2g<i<0 \\ p\mid i}} \mathscr{V}/i\mathscr{V}. \tag{3.8}$$

(Note that this last isomorphism is non-canonical: it depends on the choice of $t$.) The map $\varphi$ is the composition

$$\mathrm{H}^0(C, \Omega^1_{C/\mathscr{V}}(2g\mathscr{D})) \longrightarrow t^{-2g}\mathscr{V}[\![t]\!]dt \longrightarrow \bigoplus_{i=-2g}^{-2} \big(\mathscr{V}/(i+1)\mathscr{V}\big)\, t^i dt. \tag{3.9}$$

**3.10. Lemma.** *Consider $\Upsilon$ and $\varphi$ introduced in Notation 3.5. There is a short exact sequence*

$$0 \longrightarrow \operatorname{Ker}\varphi \longrightarrow \mathrm{H}^0(C, \Omega^1_{C/\mathscr{V}}(2g\mathscr{D})) \xrightarrow{\varphi} \Upsilon \longrightarrow 0.$$

*Let* $\operatorname{Im} d$ *be the image of the map*

$$\mathrm{H}^0(C, \mathscr{O}_C((2g-1)\mathscr{D})) \xrightarrow{\mathrm{H}^0(C,d)} \mathrm{H}^0(C, \Omega^1_{C/\mathscr{V}}(2g\mathscr{D})).$$

*Then there is an isomorphism*

$$\operatorname{Ker}\varphi \,/\, (\operatorname{Im} d \cap \operatorname{Ker}\varphi) \xrightarrow{\sim} \mathrm{H}^1_{\mathrm{dR}}(C/\mathscr{V}),$$

*which is natural in the following ways:*

*(i) The isomorphism is equivariant for maps induced by automorphisms of* $C$ *that map the divisor* $\mathscr{D}$ *to itself.*

*(ii) Let* $C' \subset C$ *be an open affine subscheme that is disjoint from the support of* $\mathscr{D}$*. There is a canonical map*

$$\mathrm{H}^0(C, \Omega^1_{C/\mathscr{V}}(2g\mathscr{D})) \longrightarrow \mathrm{H}^1_{\mathrm{dR}}(C'/\mathscr{V})$$

*which is defined as follows: one first restricts the differentials to* $C'$*, and then takes the canonical projection to de Rham cohomology (recall that on an affine scheme, de Rham cohomology is given by the closed differential forms modulo the exact ones). Restriction also gives a map*

$$\mathrm{H}^1_{\mathrm{dR}}(C/\mathscr{V}) \longrightarrow \mathrm{H}^1_{\mathrm{dR}}(C'/\mathscr{V}).$$

*These sit in a diagram*

$$\begin{array}{ccccc} \operatorname{Ker}\varphi & \twoheadrightarrow & \operatorname{Ker}\varphi/(\operatorname{Im} d \cap \operatorname{Ker}\varphi) & \xrightarrow{\sim} & \mathrm{H}^1_{\mathrm{dR}}(C/\mathscr{V}) \\ \downarrow & & & & \downarrow \\ \mathrm{H}^0(C, \Omega^1_{C/\mathscr{V}}(2g\mathscr{D})) & & \longrightarrow & & \mathrm{H}^1_{\mathrm{dR}}(C'/\mathscr{V}) \end{array}$$

*which is commutative.*

PROOF. We have already encountered the de Rham complex $\Omega^\bullet_{C/\mathscr{V}}$ in (2.1). Let us introduce two more complexes. First of all

$$\tilde{\Omega}^\bullet_{C/\mathscr{V}} : \quad \mathscr{O}_C((2g-1)\mathscr{D}) \xrightarrow{d} \Omega^1_{C/\mathscr{V}}(2g\mathscr{D}), \tag{3.11}$$

which vanishes outside degrees 0 and 1. The map $d$ is the one appearing in (3.3). There is an inclusion $\Omega^\bullet_{C/\mathscr{V}} \hookrightarrow \tilde{\Omega}^\bullet_{C/\mathscr{V}}$. The second complex we need is the cokernel complex $Q^\bullet$ of this inclusion.

From the 'stupid' short exact sequence

$$0 \longrightarrow \Omega^1_{C/\mathcal{V}}(2g\mathcal{D})[1] \longrightarrow \tilde{\Omega}^\bullet_{C/\mathcal{V}} \longrightarrow \mathcal{O}_C((2g-1)\mathcal{D})[0] \longrightarrow 0$$

one obtains an exact sequence

$$\begin{array}{l} \mathrm{H}^0(C,\mathcal{O}_C((2g-1)\mathcal{D})) \to \mathrm{H}^0(C,\Omega^1_{C/\mathcal{V}}(2g\mathcal{D})) \longrightarrow \mathrm{H}^1(C,\tilde{\Omega}^\bullet_{C/\mathcal{V}}) \longrightarrow \mathrm{H}^1(C,\mathcal{O}_C((2g-1)\mathcal{D})) \\ \to \mathrm{H}^1(C,\Omega^1_{C/\mathcal{V}}(2g\mathcal{D})) \longrightarrow \mathrm{H}^2(C,\tilde{\Omega}^\bullet_{C/\mathcal{V}}) \longrightarrow 0. \end{array}$$

Now both the spaces $\mathrm{H}^1(C,\mathcal{O}_C((2g-1)\mathcal{D}))$ and $\mathrm{H}^1(C,\Omega^1_{C/\mathcal{V}}(2g\mathcal{D}))$ vanish. Indeed, by Lemma 3.1 and flat base change, it suffices to prove this over $K$; and over $K$ it follows from Serre duality. It follows that $\mathrm{H}^2(C,\tilde{\Omega}^\bullet_{C/\mathcal{V}})=0$.

Now let $t$ be as before: a generator of $\mathcal{O}_C(-\mathcal{D})$ in a neighbourhood of $\operatorname{Supp}\mathcal{D}$. Similar to the concrete description of $\Upsilon$ above, we have that $Q^\bullet$ is the complex of skyscraper sheaves

$$\frac{t^{-(2g-1)}\mathcal{V}[\![t]\!]}{\mathcal{V}[\![t]\!]} \xrightarrow{d} \frac{t^{-2g}\mathcal{V}[\![t]\!]dt}{\mathcal{V}[\![t]\!]dt}.$$

Hence

$$\mathrm{H}^0(C,Q^\bullet) = \operatorname{Ker} d = 0$$

and

$$\mathrm{H}^1(C,Q^\bullet) = \operatorname{Coker} d = \Upsilon \oplus \mathcal{V}\frac{dt}{t}.$$

So there is an exact sequence

$$0 \longrightarrow \mathrm{H}^1_{\mathrm{dR}}(C/\mathcal{V}) \longrightarrow \mathrm{H}^1(C,\tilde{\Omega}^1_{C/\mathcal{V}}) \longrightarrow \Upsilon \oplus \mathcal{V} \longrightarrow \mathrm{H}^2_{\mathrm{dR}}(C/\mathcal{V}) \longrightarrow 0 \tag{3.12}$$

induced by $\Omega^\bullet_{C/\mathcal{D}} \rightarrowtail \tilde{\Omega}^\bullet_{C/\mathcal{D}} \twoheadrightarrow Q^\bullet$. Now $\mathrm{H}^2_{\mathrm{dR}}(C/\mathcal{V}) = \mathcal{V}$, which follows from the fact that the one dimensional free $\mathcal{V}$-module $\mathrm{H}^1(C,\Omega^1_{C/\mathcal{V}})$ maps surjectively to $\mathrm{H}^2_{\mathrm{dR}}(C/V)$ and $\mathrm{H}^2_{\mathrm{dR}}(C_K/K) \simeq K$. So the last non-zero map in this exact sequence splits. Putting everything together, we obtain a commutative diagram

$$\begin{array}{ccccccccc}
 & & & & \mathrm{H}^0(C,\mathcal{O}_C((2g-1)\mathcal{D})) & & & & \\
 & & & & \downarrow & & & & \\
 & & & & \mathrm{H}^0(C,\Omega^1_{C/\mathcal{V}}(2g\mathcal{D})) & & & & \\
 & & & & \downarrow & \searrow^{\varphi} & & & \\
0 & \longrightarrow & \mathrm{H}^1_{\mathrm{dR}}(C/\mathcal{V}) & \longrightarrow & \mathrm{H}^1(C,\tilde{\Omega}^\bullet_{C/\mathcal{V}}) & \longrightarrow & \Upsilon & \longrightarrow & 0 \\
 & & & & \downarrow & & & & \\
 & & & & 0 & & & &
\end{array}$$

in which the row and column are exact, and where the diagonal arrow is the map $\varphi$. This proves the main part of the lemma. The properties (i) and (ii) follow from the naturality of the above constructions. □

## 4. Corollaries for hyperelliptic curve

Now we will suppose that $C$ is a hyperelliptic curve over $\mathcal{V}$ of genus $g$ with (fixed) hyperelliptic involution $\iota$. Suppose the divisor $\mathcal{D}$ corresponds to a rational Weierstrass point. We will suppose that $\iota$ fixes this Weierstrass point (a vacuous condition if $g \geq 2$). Let $C' \subset C$ be an open subscheme such that $\iota$ is the identity on the complement.

As the Weierstrass points are fixed under $\iota$, there is an induced involution on all the spaces appearing in Lemma 3.10. By abuse of notation, we will denote all these induced maps again by $\iota$. If $V$ is such a space, denote by $V^{\pm}$ the subspaces $\{w \in V \mid \iota w = \pm w\}$.

**4.1. Proposition.** *The restriction map*

$$\mathrm{H}^1_{\mathrm{dR}}(C_K/K) \longrightarrow \mathrm{H}^1_{\mathrm{dR}}(C'_K/K)$$

*is injective and its image is the* $(-1)$*-eigen space for the involution induced by* $\iota$*.*

Proof. We may pass to the algebraic closure of $K$. The proposition then follows form the Gysin sequence for de Rham cohomology (see [Ha70, 8.3]). More elementarily, we may suppose we are working over the complex numbers $\mathbb{C}$. By Grothendieck's comparison theorem [Gr], it suffices to prove the analogous statement for the map

$$\mathrm{H}^1(C(\mathbb{C}),\mathbb{Q}) \longrightarrow \mathrm{H}^1(C'(\mathbb{C}),\mathbb{Q})$$

in singular cohomology. The inclusion $C'(\mathbb{C}) \subset C(\mathbb{C})$ with complement $Z$ induces a long exact sequence, part of which is

$$\mathrm{H}^1_Z(C(\mathbb{C}),\mathbb{Q}) \longrightarrow \mathrm{H}^1(C(\mathbb{C}),\mathbb{Q}) \longrightarrow \mathrm{H}^1(C'(\mathbb{C}),\mathbb{Q}) \longrightarrow \mathrm{H}^2_Z(C(\mathbb{C}),\mathbb{Q}).$$

This sequence is functorial, which implies that it preserves $(-1)$-eigenspaces for the automorphisms $\iota$. The space $\mathrm{H}^i_Z(C(\mathbb{C}),\mathbb{Q})$ is a finite sum of copies of $\mathrm{H}^i_{\{*\}}(\mathbb{R}^2,\mathbb{Q})$ with $\iota$ acting on $\mathbb{R}^2$ as an orientation preserving map. Now the last space sits in the long exact sequence of a pair

$$\mathrm{H}^0(S^1,\mathbb{Q}) \longrightarrow \mathrm{H}^1_{\{*\}}(\mathbb{R}^2,\mathbb{Q}) \longrightarrow \mathrm{H}^1(\mathbb{R}^2,\mathbb{Q}) \longrightarrow \mathrm{H}^1(S^1,\mathbb{Q}) \longrightarrow \mathrm{H}^2_{\{*\}}(\mathbb{R}^2,\mathbb{Q}) \longrightarrow \mathrm{H}^2(\mathbb{R}^2,\mathbb{Q})$$

and hence its $(-1)$-eigenspace vanishes for $i = 1,2$ as it is wedged in between zero-spaces. Therefore, $\mathrm{H}^1(C(\mathbb{C}),\mathbb{Q})^- \xrightarrow{\sim} \mathrm{H}^1(C'(\mathbb{C}),\mathbb{Q})^-$. Over an algebraically closed field, $C$ modulo the action of $\iota$ is isomorphic to $\mathbb{P}^1$. Taking cohomology with $\mathbb{Q}$ coefficients commutes with taking invariants under a finite group action. Hence $\mathrm{H}^1(C(\mathbb{C}),\mathbb{Q})^+ = \mathrm{H}^1(C(\mathbb{C}),\mathbb{Q})^{\langle\iota\rangle} = \mathrm{H}^1(\mathbb{P}^1(\mathbb{C}),\mathbb{Q}) = 0$. □

**4.2. Theorem.** *The restriction to* $(\operatorname{Ker}\varphi)^-$ *of the map of Lemma 3.10 is injective:*

$$(\operatorname{Ker}\varphi)^- \hookrightarrow \mathrm{H}^1_{\mathrm{dR}}(C/\mathcal{V})$$

*and its image contains* $2\mathrm{H}^1_{\mathrm{dR}}(C/\mathcal{V})$. *In particular:*

*(i) there is an isomorphism*

$$\mathrm{H}^0(C_K,\Omega^1_{C_K/K}(2g\mathcal{D}_K))^- \xrightarrow{\sim} \mathrm{H}^1_{\mathrm{dR}}(C_K/K)$$

*compatible with the canonical maps to* $\mathrm{H}^1_{\mathrm{dR}}(C'_K/K)$;

*(ii) if* $p \neq 2$, *then* $(\operatorname{Ker}\varphi)^- \simeq \mathrm{H}^1_{\mathrm{dR}}(C/\mathcal{V})$ *(recall that* $p$ *denotes the residue characteristic of* $\mathcal{V}$*).*

Proof. By the gap sequence for hyperelliptic curves (see section 1), the involution $\iota$ acts as the identity on $\mathrm{H}^0(C_K,\mathcal{O}_{C_K}((2g-1)\mathcal{D}_K))$. So by flat base change and Lemma 3.1, $\iota$ acts also as the identity on $\mathrm{H}^0(C,\mathcal{O}_C((2g-1)\mathcal{D}))$. This space surjects onto the kernel of the map $(\operatorname{Ker}\varphi) \to \mathrm{H}^1_{\mathrm{dR}}(C/\mathcal{V})$ of Lemma 3.10, so $(\operatorname{Ker}\varphi)^- \to \mathrm{H}^1_{\mathrm{dR}}(C/\mathcal{V})$ is injective.

By Proposition 4.1 and Proposition 3.2, $\iota$ acts as $-1$ on $\mathrm{H}^1_{\mathrm{dR}}(C/\mathcal{V})$. If $v \in \mathrm{H}^1_{\mathrm{dR}}(C/\mathcal{V})$, it lifts to an element $\tilde{v}$ of $\operatorname{Ker}\varphi$ by Lemma 3.10. But $2\tilde{v} = (\tilde{v}+\iota\tilde{v}) + (\tilde{v}-\iota\tilde{v})$. So $\tilde{v}-\iota\tilde{v} \in (\operatorname{Ker}\varphi)^-$ maps to $2v$. □

## 5. Application to Kedlaya's algorithm ($p$ odd)

We are now ready to apply the above statements to Kedlaya's algorithm. We will not describe Kedlaya's algorithm in detail (see [Ked] or [Edix]), but concentrate on the parts that we want to modify. Throughout this section, we assume $p \neq 2$: so the residue field $k$ of $\mathcal{V}$ has odd characteristic.

The input of the algorithm is a hyperelliptic curve over $k$ of genus $g$ with a rational Weierstrass point. It is lifted to a hyperelliptic curve $C$ over $\mathcal{V}$ with a divisor $\mathcal{D}$ (the 'section at infinity'). The complement $C^0$ of the support of $\mathcal{D}$ in $C$ is described as follows: it is the closed subscheme of $\mathbb{A}^2_{\mathcal{V}}$ given by the equation

$$y^2 = Q(x), \tag{5.1}$$

where $Q(x) \in \mathcal{V}[x]$ is a monic polynomial of degree $2g+1$ whose reduction mod $p$ has no double roots in an algebraic closure of $k$. This last property assures that $C^0$ is smooth over $\mathcal{V}$. The hyperelliptic involution is given by $x \mapsto x$, $y \mapsto -y$. One denotes by $C' = D(y)$ the complement of $x$-axis in $C^0$. Note that $\iota$ acts as the identity on $C \setminus C'$.

**5.2. Proposition.** *The elements*

$$\frac{dx}{y},\ x\frac{dx}{y},\ \ldots,\ x^{2g-1}\frac{dx}{y} \tag{5.3}$$

*of* $\Omega^1_{K(C)/K}$ *freely generate the* $\mathcal{V}$*-submodule* $\mathrm{H}^0(C,\Omega^1_{C/\mathcal{V}}(2g\mathcal{D}))^-$. *The order at* $\mathcal{D}$ *of such an element* $x^i\frac{dx}{y}$ *is equal to* $2g-2-2i$.

PROOF. The module is free by Lemma 3.1. Let $A = \mathcal{V}[x,y]/(y^2 - Q(x))$ be the coordinate ring of $C^0$. Then

$$A^+ = \bigoplus_{i\geq 0} \mathcal{V}x^i \qquad \text{and} \qquad \Omega_{A/\mathcal{V}} \simeq (A dx \oplus A dy)/(2y dy - Q'(x)dx).$$

By the Jacobi criterium, $D(y)$ and $D(Q'(x))$ cover $C^0 = \operatorname{Spec} A$. On the first subspace $dx/y$ is a generator, while on the second $2dy/Q'(x)$ is, and we can glue these to a global section of $\Omega_{A/\mathcal{V}}$. As $\iota(dx/y) = -dx/y$, it follows that $\mathrm{H}^0(C^0,\Omega^1_{C/\mathcal{V}})^- = \bigoplus_{i\geq 0} \mathcal{V}x^i dx/y$. It remains to see which part of this space extends to a section of $\mathrm{H}^0(C,\Omega^1_{C/\mathcal{V}}(2g\mathcal{D}))$. The complement $C \setminus C^0$ is the support of $\mathcal{D}$ and consists of a generic point $\xi$ and a special point. An order calculation gives

$$\operatorname{ord}_\xi(x^i dx/y) = 2g-2-\sum_{Q\in C^0(\overline{K})} \operatorname{ord}_Q(x^i dx/y) = 2g-2-2i$$

(where $K \subset \overline{K}$ is an algebraic closure). By normality, a local section of $\Omega^1_{C/\mathcal{V}}(2g\mathcal{D})$ that extends to $\xi$ also extends to the special point. So the submodule of $\mathrm{H}^0(C^0,\Omega^1_{C/\mathcal{V}})^-$ of sections that extend to global sections of $\Omega^1_{C/\mathcal{V}}(2g\mathcal{D})$ is indeed spanned by (5.3). □

By Theorem 4.2 the elements (5.3) form a basis for the $K$-vector space $\mathrm{H}^1_{\mathrm{dR}}(C_K/K)$, and so by Proposition 4.1 they also form a basis for $\mathrm{H}^1_{\mathrm{dR}}(C'_K/K)^-$. (By the way, this last result can also be obtained by explicitly calculating the cokernel of the derivation

$$A'_K \xrightarrow{d} \Omega_{A'_K/K},$$

where $A'_K = K[x,y,y^{-1}]/(y^2 - Q(x))$ is the coordinate ring of the affine space $C'_K$.) As we have seen in Theorem 2.6, Frobenius acts on $\mathrm{H}^1_{\mathrm{dR}}(C'_K/K)^- = \mathrm{H}^1_{\mathrm{dR}}(C_K/K)$ in a natural way.

Taking up Kedlaya's algorithm again, it now calculates a matrix $M$ which represents the automorphism induced by Frobenius with respect to the basis (5.3). The coefficient of $M$ are in $K$. (In the actual algorithm, $M$ is calculated with a certain $p$-adic precision. We can ignore this fact, as it is irrelevant for our present purposes.) As a first application of Lemma 3.10, we deduce the following fact about this matrix.

**5.4. Proposition.** *Put* $r = \lfloor \log_p(2g-1) \rfloor$. *The coefficients of the matrix* $p^r M$ *are all elements of* $\mathcal{V}$.

PROOF. Denote by $\Gamma$ the automorphism induced by Frobenius on $\mathrm{H}^1_{\mathrm{dR}}(C_K/K)$ and, by slight abuse of notation, on its subspace $\mathrm{H}^1_{\mathrm{dR}}(C/\mathcal{V})$. There are the inclusions

$$(\mathrm{Ker}\,\varphi)^- \hookrightarrow \mathrm{H}^0(C,\Omega^1_{C/\mathcal{V}}(2g\mathcal{D}))^- \hookrightarrow \mathrm{H}^0(C_K,\Omega^1_{C_K/K}(2g\mathcal{D}_K))^-.$$

Theorem 4.2 says that the space on the left is canonically isomorphic to $\mathrm{H}^1_{\mathrm{dR}}(C/\mathcal{V})$ and that the space on the right is canonically isomorphic to $\mathrm{H}^1_{\mathrm{dR}}(C_K/K)$; hence $\Gamma$ acts on the two spaces at the left and the right.

Let $a$ be an integer satisfying $0 \le a \le 2g-1$. Denote by $\lambda_i \in K$ (with $i = 0, \ldots, 2g-1$) the matrix elements for which $\Gamma(x^a \frac{dx}{y}) = \sum_{i=0}^{2g-1} \lambda_i x^i \frac{dx}{y}$. From semi-linearity it follows that $\Gamma(p^r x^a \frac{dx}{y}) = \sum_{i=0}^{2g-1} p^r \lambda_i x^i \frac{dx}{y}$. But Lemma 3.10 says that $p^r x^a \frac{dx}{y}$ belongs to the kernel of $\varphi$. Hence also $\Gamma(p^r x^a \frac{dx}{y})$ belongs to this kernel. We conclude from Proposition 5.2 that $p^r \lambda_i \in \mathcal{V}$ for all $i$. □

We will now describe how to determine elements

$$e_j = \sum_{i=0}^{2g-1} c_{ij} x^i \frac{dx}{y} \qquad (j = 0,\ldots,2g-1, \quad c_{i,j} \in \mathcal{V})$$

that form a $\mathcal{V}$-basis for $\mathrm{H}^1_{\mathrm{dR}}(C/\mathcal{V})$, or, what by Theorem 4.2 amounts to the same thing, that form a basis for the free submodule $(\operatorname{Ker}\varphi)^-$ of $\mathrm{H}^0(C,\Omega^1_{C/\mathcal{V}}(2g\mathscr{D}))^-$. The Frobenius with respect to this basis will then have coefficients in $\mathcal{V}$.

First we consider an equation for $C$ in a neighbourhood of $\mathscr{D}$: choosing coordinates $u = 1/x$ and $t = yu^{g+1}$, (5.1) gives

$$t^2 = P(u) := Q(1/u)u^{2g+2}.$$

As $Q(x)$ is monic of degree $2g+1$, we can write $P(u) = uR(u)$, with $R(u) \in \mathcal{V}[u]$ satisfying $R(0) = 1$. So also $P'(0) = 1$ (and $P$ mod $p$ has no double roots). The divisor $\mathscr{D}$ corresponds to $u = 0$, $t = 0$.

Recall that the map $\varphi$ is the map (see (3.9)) given by truncation

$$\mathrm{H}^0(C,\Omega^1_{C/\mathcal{V}}(2g\mathscr{D})) \longrightarrow \bigoplus_{i=-2g}^{-2} \left(\mathcal{V}/(i+1)\mathcal{V}\right) t^i dt.$$

Let $\varphi^-$ be the restriction of $\varphi$ to the $(-1)$-eigenspace for the action of the involution $\iota$. Since $\iota(t) = -t$, we can consider $\varphi^-$ as a map

$$\mathrm{H}^0(C,\Omega^1_{C/\mathcal{V}}(2g\mathscr{D}))^- \longrightarrow \bigoplus_{\lambda=1}^{g} \left(\mathcal{V}/(2\lambda-1)\mathcal{V}\right) t^{-2\lambda} dt,$$

where we have conveniently relabeled the summation. For future use, let us also introduce the canonical projections

$$\bigoplus_{\lambda=1}^{g} \left(\mathcal{V}/(2\lambda-1)\mathcal{V}\right) t^{-2\lambda} dt \xrightarrow{\mathrm{pr}_\varkappa} \left(\mathcal{V}/(2\varkappa-1)\mathcal{V}\right).$$

where $\varkappa = 1, \ldots, g$.

Considering the order at $\mathscr{D}$ of the elements $\frac{dx}{y}, \ldots, x^{g-1}\frac{dx}{y}$ given by Proposition 5.2, we observe that these elements are mapped by $\varphi$ to zero. The image of the other elements in (5.3) can be calculated as follows. Set $F = \mathcal{V}/p^r\mathcal{V}$, with $r = \lfloor \log_p(2g-1) \rfloor$. Determine $R(u)^{-1}$ as an element of $F[[u]]/(u^g)$, which is possible since $R(0) = 1$. By iterating the relation

$$u = t^2 R(u)^{-1} \qquad \left[ \, = t^2(1 + \mathcal{O}(u)) \, \right]$$

we obtain an expression for $u$ as element of $F[[t]]/(t^{2g-1})$. The relation

$$x^{g-1}\frac{dx}{y} = -2P'(u)^{-1} dt$$

enables us to express $x^{g-1}\frac{dx}{y}$ as an element of $F[\![t]\!]dt/(t^{2g-1})$; it even sits in $(-2dt + tF[\![t]\!]dt)/(t^{2g-1})$ and note that $-2$ is a unit in $F$. Inductively for $\varkappa = 1, \ldots, g$, we can determine the image of $x^{\varkappa+g-1}\frac{dx}{y}$ in $t^{-2\varkappa}F[\![t]\!]dt/(t^{2g-1-2\varkappa})$ by multiplying $x^{\varkappa+g-2}\frac{dx}{y}$ by $x = t^{-2}R(u)$. The canonical projection

$$t^{-2g}F[\![t]\!]dt/(t^{-1}) \longrightarrow \bigoplus_{\lambda=1}^{g} \left(\mathcal{V}/(2\lambda-1)\mathcal{V}\right) t^{-2\lambda}dt$$

then gives an expression for $\varphi^{-}(x^i\frac{dx}{y})$ with $g \le i < 2g$. Observe that $\mathrm{pr}_{i-g+1}\varphi(x^i\frac{dx}{y})$ is a unit in $\left(\mathcal{V}/(2(i-g+1)-1)\mathcal{V}\right)$ (where 0 is a unit if this space is trivial) and that $\mathrm{pr}_{\varkappa}\varphi(x^i\frac{dx}{y}) = 0$ for $\varkappa > i-g+1$.

For $i = 0, \ldots, 2g-1$, put $V_i = \mathcal{V}\frac{dx}{y} \oplus \mathcal{V}x\frac{dx}{y} \oplus \cdots \oplus \mathcal{V}x^i\frac{dx}{y}$. This defines a filtration on $V_{2g-1} = \mathrm{H}^0(C, \Omega^1_{C/\mathcal{V}}(2g\mathcal{D}))^{-}$. Put $W_i = (\mathrm{Ker}\,\varphi^{-}) \cap V_i$, which defines the corresponding filtration on $(\mathrm{Ker}\,\varphi)^{-} = (\mathrm{Ker}\,\varphi^{-})$. We will define a basis $e_0, \ldots, e_i$ of $W_i$ inductively.

For $0 \le j \le g-1$, put $e_j = x^j\frac{dx}{y}$. As we already observed that $\varphi(e_j) = 0$, these elements provide bases for the $W_i$ with $i \le g-1$. Now consider $i > g-1$ and suppose that we have defined elements $e_0, \ldots, e_{i-1}$ that form a basis of $W_{i-1}$. Let $\varkappa$ be such that $i = \varkappa + g - 1$ and consider the element $(2\varkappa-1)x^i\frac{dx}{y}$. For $\lambda \ge \varkappa$ we then have $\mathrm{pr}_{\lambda}\varphi^{-}((2\varkappa-1)x^i\frac{dx}{y}) = 0$. Since $\varphi^{-}(V_{i-1}) = \bigoplus_{\lambda=1}^{\varkappa-1}\left(\mathcal{V}/(2\lambda-1)\mathcal{V}\right) t^{-2\lambda}dt$, we can find elements $c_j \in \mathcal{V}$ such that $e_i = (2\varkappa-1)x^i\frac{dx}{y} + \sum_{j=0}^{i-1} c_j x^j\frac{dx}{y}$ sits in the kernel of $\varphi^{-}$. So $e_i \in W_i$.

Consider an element $x = \sum_{j=0}^{i} d_j x^j\frac{dx}{y} \in W_i$. As $0 = \mathrm{pr}_{\varkappa}\varphi^{-}(x) = \mathrm{pr}_{\varkappa}\varphi^{-}(d_i x^i\frac{dx}{y}) = d_i\mathrm{pr}_{\varkappa}\varphi^{-}(x^i\frac{dx}{y})$ and $\mathrm{pr}_{\varkappa}\varphi^{-}(x^i\frac{dx}{y})$ is a unit, it follows that $d_i = (2\varkappa-1)d'$ for some element $d' \in \mathcal{V}$. Now $x - d'e_i \in W_{i-1}$, so by the induction hypotheses $x \in W_{i-1} \oplus \mathcal{V}e_i$. Thus the elements $e_0, \ldots, e_i$ form a basis of $W_i$.

## 6. Comments on the case of characteristic 2

Denef and Vercauteren [D-V] have adapted Kedlaya's algorithm to the case of characteristic 2. Let us briefly sketch what our methods can accomplish in this case.

From now on we assume $p = 2$. An equation for $C$ outside $\mathcal{D}$ is given by

$$y^2 + h(x)y = f(x),$$

where $f(x) \in \mathcal{V}[x]$ is monic of degree $2g+1$ and $h(x) \in \mathcal{V}[x]$ has degree $\le g$ (in fact, one can assume that $h(x)$ satisfies more conditions; see [ibid.]). As in the proof of Proposition 5.2, one shows that the elements

$$\frac{dx}{2y+h(x)},\ x\frac{dx}{2y+h(x)},\ \ldots,\ x^{2g-1}\frac{dx}{2y+h(x)} \tag{6.1}$$

form a basis for the free $\mathscr{V}$-submodule $\mathrm{H}^0(C,\Omega^1_{C/\mathscr{V}}(2g\mathscr{D}))^-$. A method similar to the algorithm in the previous section gives a change of basis from (6.1) to one that spans $(\mathrm{Ker}\,\varphi)^-$. But unlike before, we cannot conclude that this basis generates $\mathrm{H}^1_{\mathrm{dR}}(C/\mathscr{V})$. All we can conclude from Theorem 4.2 is:

**6.2. Proposition.** *Let $M'$ be the matrix of Frobenius with respect to the basis for* $\mathrm{Ker}\,\varphi$. *Then the coefficients of $2M'$ are in $\mathscr{V}$.*

There is one further thing to do: the basis used in [ibid.] for $\mathrm{H}^1_{\mathrm{dR}}(C'_K/K)$ is not the one given by (6.1), so one has to work out the transition matrix by a calculation in affine de Rham cohomology.

### References for chapter 3